Volume 96
Number 2
Summer 2010# Journal of the WASHINGTON ACADEMY OF SCIENCES

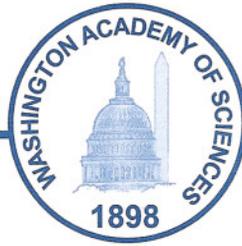

Editor's Comments  *J. Maffucci* ............................................................................. i
Letters to the Editor ............................................................................................ iii
"Human Music" a Theoretical Model of How Music Induces Affect  *D. Teie* ............... 1
A New Perspective on the Early History of the American Society for Cybernetics
    *E. Corona and B. Thomas* ........................................................................ 21
Air Traffic Controller Workload: Estimating Look-Ahead Conflict Detection Counts
    *N. Coleman and E. Feldman* ................................................................... 35
A Digital-Discrete Method For Smooth-Continuous Data Reconstruction  *L. Chen* ......... 47
Outgoing President's Speech  *Kiki Ikossi* ................................................................ 67
Incoming President's Speech  *Mark Holland* ........................................................... 72
Banquet 2010 photos ......................................................................................... 73
ISSN 0043-0439       Issued Quarterly at Washington DC

# Digital-Discrete Method For Smooth-Continuous Data Reconstruction

Li Chen
Department of Computer Science and Information Technology
University of the District of Columbia

## Abstract

A systematic digital-discrete method for obtaining continuous functions with smoothness to a certain order ($C^{(n)}$) from sample data is designed. This method is based on gradually varied functions and the classical finite difference method. This new method has been applied to real groundwater data and the results have validated the method. This method is independent from existing popular methods such as the cubic spline method and the finite element method. The new digital-discrete method has considerable advantages for a large number of real data applications. This digital method also differs from other classical discrete methods that usually use triangulations. This method can potentially be used to obtain smooth functions such as polynomials through its derivatives $f^{(k)}$ and the solution for partial differential equations such as harmonic and other important equations.

## Introduction

OF THE MOST COMMON PROBLEMS in data reconstruction is to fit a tion based on the observations of some sample (guiding) points. The l 'some' is important here. When one has all the data points in hand, fitting a function to them is trivial. When the data set is incomplete, one wishes to reconstruct the data, then building a proper fitting ion is more difficult. In this paper, we present recently developed ithms for smooth-continuous data reconstruction based on the digital- ete method. The classical discrete method for data reconstruction is l on domain decomposition according to guiding (or sample) points, hen the Spline method (for polynomial fitting) or finite elements od (for Partial Differential Equations) is used to fit the data. Some ssful methods have been discovered or proposed to solve the em including the Voronoi-based surface method[1] and the moving square method [2][12][15][18][19]. A comprehensive review was nted in [4].

---

thematics, a Voronoi diagram is a special kind of decomposition of a metric space ined by distances to a specified discrete set of objects in the space, *e.g.,* by a e set of points. It is named after Georgy Voronoi (1868-1908).

# A Digital-Discrete Method For Smooth-Continuous Data Reconstruction

Li Chen
Department of Computer Science and Information Technology
University of the District of Columbia
lchen@udc.edu

**Abstract**

*A systematic digital-discrete method for obtaining continuous functions with smoothness to a certain order ($C^{(n)}$) from sample data is designed. This method is based on gradually varied functions and the classical finite difference method. This new method has been applied to real groundwater data and the results have validated the method. This method is independent from existing popular methods such as the cubic spline method and the finite element method. The new digital-discrete method has considerable advantages for a large number of real data applications. This digital method also differs from other classical discrete methods that usually use triangulations. This method can potentially be used to obtain smooth functions such as polynomials through its derivatives $f^{(k)}$ and the solution for partial differential equations such as harmonic and other important equations.*

## 1. Introduction

One of the most common problems in data reconstruction is to fit a function based on the observations of some sample (guiding) points. The word 'some' is important here. When one has all the data points in hand, then fitting a function to them is trivial. When the data set is incomplete, and one wishes to reconstruct the data, then building a proper fitting function is more difficult. In this paper, we present recently developed algorithms for smooth-continuous data reconstruction based on the digital-discrete method. The classical discrete method for data reconstruction is based on domain decomposition according to guiding (or sample) points, and then the Spline method (for polynomial fitting) or finite elements method (for Partial Differential Equations) is used to fit the data. Some successful methods have been discovered or proposed to solve the problem including the Voronoi-based surface method[1] and the moving least square method [2][12][15][18][19]. A comprehensive review was presented in [4].

Our method is based on the *gradually varied function* (see the appendix for a definition) that does not assume the property of the linear separability among guiding points, i.e. no domain decomposition methods are needed [8][9]. We design a systematic digital-discrete method for obtaining continuous functions with smoothness to a certain order $(C^{(n)})$[2] from sample data. This design is based on gradually varied functions and the classical finite difference method. We also demonstrate the flexibility of the new method and its potential to solve a variety of problems. The examples include some real data from water well logs and harmonic functions on closed 2D manifolds. These validate the method. We present several different algorithms. This method can be easily extended to higher multi-dimensions.

---

[1] In mathematics, a Voronoi diagram is a special kind of decomposition of a metric space determined by distances to a specified discrete set of objects in the space, *e.g.,* by a discrete set of points. It is named after Georgy Voronoi.

[2] $C^0$ is the class of all continuous functions. $C^k$ is the class of differentiable functions whose kth derivative is continuous.

This method is independent from existing popular methods such as the cubic spline method and the finite element method. The new digital-discrete method has considerable advantages for a large number of real data applications. This digital method also differs from other classical discrete methods that usually use triangulations. This method can potentially be used to obtain smooth functions such as polynomials through its derivatives $f^{(k)}$ and the solution for partial differential equations such as harmonic and other important equations [3][4][5].

## 2. Background and Basic Concepts

We will deal with the following two real world problems: (1) Given a set of points and its observation (function) values at these points, extend the values to a larger set. (2) When observing an image, if an object is extracted from the image, a representation of the object can sometimes be described by its boundary curve. If all the values on the boundary are the same, then we can restore the object by filling the region. If the values on the boundary are not the same and if we assume the values are "continuous" on the boundary, then we need a fitting algorithm to find a surface. Both problems involve extending the original dataset to the entire region. We will use two real data examples in the following sections to explain these.

We address the following question about extending the data: Let $D$ be a domain and $J$ be a subset of $D$. If $f$ is "continuous" or "smooth" on $J$, is there a general method that yields an extension $F$ of $f$ for set $D$ that has the continuity or smoothness property?

In continuous space, this problem is related to the Dirichlet Problem (when $J$ is the boundary) and the Whitney extension problem (when $J$ is a subspace of $D$, the space $R^{(n)}$) [13][14][19].

Why are the existing numerical methods not perfect? Here is an explanation using splines. We show in Fig 2.1 an example that contains four sample data points. If the boundary were irregular, we would need to use a 2D B-Spline to divide the boundary into four segments. The different partitions would yield different results (One is free to do this in different ways with different choices). Fig 2a shows one choice for a linear interpolation. Fig 2b shows another choice. Fig. 2c shows how a gradually varied interpolation will fit the data. If we have five sample points, we would have 10 different piece-wise linear interpolations. For six points, we may have more than 30 piece-wise linear functions [4]. For more literature review, see [4]. Here we just use an example to illustrate the differences that can arise with choice of partition. [4].

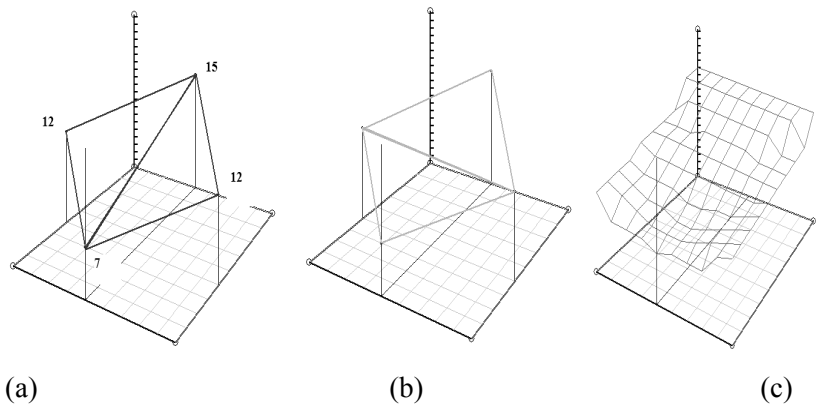

(a)     (b)     (c)

Fig.2.1. (a) (b) Two piecewise linear interpolations, we do not know which one is correct. (c) The gradually varied interpolation result shows a quite reasonable non-linear fitting.

To answer the question posed earlier: What is a "continuous" or "smoothing" function in discrete space? We have defined the so-called *gradually varied functions* (GVF) for the purpose of constructing continuous functions in discrete space [6][7][8][9][10]. See the Appendix for the mathematical definition.

The basic concept of gradual variation is to define small changes between two points in a discrete domain that can be built on any graph [6][7]. So, a gradually varied surface is a special discrete surface. In gerneal, A digital surface is formed by the moving of a line segment.

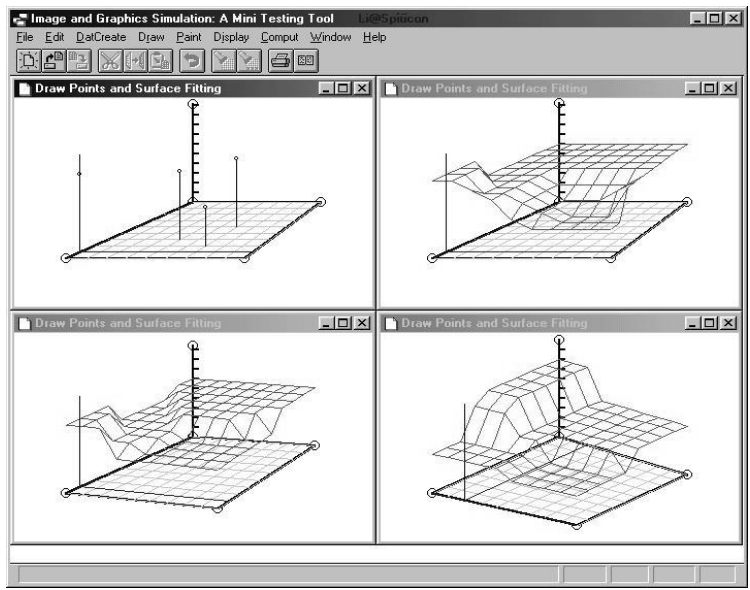

Fig. 2.2 Examples of gradually varied functions

In 1997 a gradually varied surface fitting software component was included in a lab-use software system. The software demonstrates the arbitrary guiding points of a gradually varied surface fitting in a 10x10 grid domain (Fig. 2.2).

The gradually varied function is tightly related to the Lipschitz function and the local Lipschitz function. A brief comparison of the method of gradually varied functions and the McShane-Whitney extension method is in [3]. In theory, McShane and Whitney obtained an important theorem for the Lipschitz function extension [17][22]. Kirszbraun and later Valentine studied the Lipschitz mapping extension for Hilbert spaces [21]. For more information, see [13][3].

Theorem 2.1 (see Appendix) can be used for a single surface fitting if the condition in the theorem is satisfied. A problem occurs when the sample data does *not* satisfy the condition of fitting. Then the original algorithm cannot be used directly for individual surface fitting. Another problem is that the theorem is only for "continuous" surfaces. It does not suggest a solution for differentiable or smooth functions.

On the other hand, we can define the discrete immersion problem as follows: Let $D_1$ and $D_2$ be two discrete manifolds (for instance, piecewise linear approximations of topological manifolds) and $f: D_1 \rightarrow D_2$ be a mapping. $f$ is said to be an immersion from $D_1$ to $D_2$ or a gradually varied operator if x and y are adjacent in $D_1$ implying $f(x) = f(y)$, or $f(x), f(y)$ are adjacent in $D_2$ [6][7]. We can define immersion-extendable and the normal immersion (the gradually varied extension) [6][7].

## 3. Smooth Gradually Varied Functions

The fitting algorithm for a typical gradually varied function does not really need the functions $A_1, A_2, \ldots A_n$ and only needs the real numbers $1, \ldots, n$ as the graph-theoretic solution. However, information about the rational or real numbers of $A_1, A_2, \ldots A_n$ are critical factors in the smooth gradually varied functions. In other words, in order to obtain a gradually varied interpolation or uniform approximation, for instance in a 2D grid space, we do not have to use the actual $A_i$. We only need to use $i, i=1,\ldots,n$. This is because gradual variation is an abstract concept. However, to consider smoothness, we have to consider the derivatives that will involve the actual values $A_i$, $i=1,\ldots,n$.

The key to the method for reconstructing a smooth gradually varied function is first to calculate a continuous function (gradually varied function), then to obtain the partial derivatives, and finally to modify the original function. A method to keep the partial derivative functions gradually varied is designed as the necessary part for this new method. This procedure can be recursively done in order to obtain high order derivatives. Then we can use the Taylor expansion for local fittings. The use of the Taylor expansion for 2D surfaces was designed by many researchers [16][13].

### 3.1 The Basics of the New Method
Given $J \subseteq D$, and $f_J: J \rightarrow \{A_1, A_2, \ldots A_n\}$. Let $f_D$ be a gradually varied extension of $f_J$ on $D$, which is a simply connected region in 2D grid space. We can calculate the derivatives

$$f_D'|_x = \tfrac{\partial}{\partial x} f_D(x,y) = f_D(x+1, y) - f_D(x, y) \tag{3.1}$$

$$f_D'|_y = \tfrac{\partial}{\partial y} f_D(x,y) = f_D(x, y+1) - f_D(x, y) \tag{3.2}$$

These derivatives will be regarded as an estimation of the fitted surfaces. Sometimes we know the values of the derivatives for the whole or subset of the domain. Then we will use those instead of the above equations. Sometimes, if we are not confident with the whole function (e.g. calculated by the above equations), we will use gradually varied functions to fit the samples of derivative values. After we have the derivatives, $f'$, we can then use them to re-calculate or update the original (gradually varied) fitted function by adding the first derivative component.

The above method can be used to calculate the different orders of derivatives. Then, the Taylor expansion at the sampling point will be applied to a region with a certain radius.

### 3.2 Calculation of the Derivatives using Gradually Varied Functions
One of the main objectives of this paper is to obtain "continuous" derivative functions using gradually varied reconstruction. After a function is reconstructed, we can then get all the orders of derivatives. There is usually no need to use fitting methods for obtaining derivative functions. However, in order to maintain the continuity of the derivative functions, we either need to smooth the derivative function or use another method to make a continuous function. In some situations, from the observing data point, we can get not only the values of the function itself, but also the values of the derivative function. For example, in the groundwater well log data, we can get both the water lever and the speed of the water flow. For some boundary value problems, we could get the boundary derivatives (called the Neumann problem in partial differential equations).

This method will directly use gradually varied fitting to get gradually varied $f_x$ and $f_y$; we can use the same technique to obtain $f_{xx}$, $f_{xy}$, and $f_{yy}$. So eventually we can get every order of derivatives we want.

Iteration with special treatments may be needed for a good fit. For instance, we usually need to iterate at a lower order of the derivative function until it is stabilized before calculating higher order derivatives.

This method differs fundamentally from Fefferman's theoretical method in digital reconstruction, which uses a system of linear inequalities and an objective function (this is called linear programming) to find the solution at each point. The inequalities are for all different orders of derivatives [13]. However, for just the $C^{(0)}$ function, the objective of Fefferman's method is very similar to that of the gradually varied reconstruction method.

### 3.3 Recalculation of the Function using Taylor Expansion

After different derivatives are obtained, we can use Taylor expansion to update the value of the gradually varied fitted function (at $C^{(0)}$). In fact, at any order $C^{(k)}$, we can update it using a higher order of derivatives.

The Taylor expansion is based on the formula of the Taylor series, which has the following generalized form:

$$f(x_1,\cdots,x_d) = \sum_{n_1=0}^{\infty} \cdots \sum_{n_d=0}^{\infty} \frac{(x_1-a_1)^{n_1}\cdots(x_d-a_d)^{n_d}}{n_1!\cdots n_d!}\left(\frac{\partial^{n_1+\cdots+n_d}f}{\partial x_1^{n_1}\cdots\partial x_d^{n_d}}\right)(a_1,\ldots,a_d). \quad (3.3)$$

For example, for a function that depends on two variables, $x$ and $y$, the Taylor series of the second order using the guiding point $(x_0, y_0)$ is:

$$f(x,y) \approx f(x_0,y_0) + (x-x_0)f_x(x_0,y_0) + (y-y_0)f_y(x_0,y_0)$$

There are several ways of implementing this formula. We have chosen the ratio of less than half of the change. An iteration process is designed to make the new function converge.

## 4. Algorithm Design

In Section 3, a systematic digital-discrete method for obtaining continuous functions with smoothness to a certain order ($C^n$) from the sample data was defined. In order to implement this method, we will now design the new algorithms to accomplish our task.

### 4.1 The Main Algorithm

The new algorithm tries to search for the best fit. We have added a component of the classical finite difference method in order to obtain derivatives for the smooth fitting. Start with a particular dataset consisting of guiding points defined on a grid space. The major steps of the new algorithm are as follows (This is for 2D functions. For 3D function, we would only need to add a dimension):

**Step 1:** *Load guiding points. In this step we load the data points with observation values.*
**Step 2:** *Determine the resolution. Locate the points in grid space.*
**Step 3:** *Function extension according to Theorem 2.1. This way, we obtain gradually varied or near gradually varied (continuous) functions. In this step the local Lipschitz condition is used.*
**Step 4:** *Use the finite difference method to calculate partial derivatives. Then obtain the smoothed function.*
**Step 5:** *Some multilevel and multi resolution method may be used to do the fitting when data set is large.*

## 4.2 The Iterated Method for Each Step of Calculating Derivatives

Iterate to find the fitting function – step 3 of the algorithm. After we have obtained $f_x, f_y$ then we can re-compute the gradually varied fitting functions (*GVF*) using $f_x, f_y$ to direct the final fitting output. Every time we need to update the result until the new function has no change. Then we have a fixed $f_x$ or $f_y$. So we can do a gradually varied fitting on the selected points in $f_x$ or $f_y$, before repeating. We will obtain $f_{xx}, f_{xy}$, and $f_{yy}$. Update $f_x$ (or $f_y$), until $f_x$ has no change. Then we return back to change $f$. If we knew $f_x, f_y$, we can use $f_x$ and $f_y$ to guide the fitting.

In other words, this method uses $F_x$, $F_y$, or GVF$(f_x)$ and GVF$(f_y)$ to update the $F$ = GVF$(f)$. We use numerical updates unlike the first gradually varied fitting where we use "digital" fitting ($A_i$, $i=1,\ldots,n$), iterate based on the fitting orders.

We can either -- Choice (A): Update the whole $F$, and then compute $F_x$, $F_y$. Then repeat until there is no change or very small change/error within a threshold. Or -- Choice (B): Update based on each order (with respect to the distance to the guiding points) then editing $F_x$, $F_y$ using the updated versions as guiding points.

Repeat until there is no change, and we get $F_x$, $F_y$. or GVF$(f_x*)$, GVF$(f_y*)$ (* means some extreme points are added). Then we compute $F_{xx}$, $F_{xy}$, and $F_{yy}$, and so on and so forth. Using $F_{xx}$, etc to update $F_x$ and $F_y$, and then back to updating $F$ again.

## 4.3 The Multi-scaling Method[3]

The multi scaling method is to choose a base scale and then refining the scale by 2. This is just like the wavelet method. Other multi-scaling methods for PDE can be found in [23].

If there is more than one guiding point in a pixel or block unit, we can chose one or use its average value. Computing the gradually varied function at scale $k$ gives us $F_k$=GVF$(f,k)$, we can then get the $F_{(x,k)}$. We can calculate and insert the value at 1/2 point surrounding the guiding points (do the corresponding process if the pixel contains more points, restores the points or uses the average value).

Computing the whole insertion or computing it in an order by surrounding points then uses the new calculated points as guiding points for the gradually varied functions. The resulting gradually varied function is in the new scale. This will guarantee the derivatives at guiding points.

---

[3] The Multiscale method is a class of algorithmic techniques for solving efficiently and effectively large-scale computational and optimization problems. The main objective of a multiscale algorithm is to create a hierarchy of problems (*coarsening*), each representing the original problem, but with fewer degrees of freedom.

Then we use this function to refine the scale again by keeping the inserted points as guiding points. We will now have more points surrounding the original guiding points, and so on and so forth. We can get our $F$ in a predefined scale. This will obtain a good derivative as well. Using higher order derivatives can be obtained recursively.

We can also use this to calculate $f_{xx}$, $f_{xy}$, and $f_{yy}$, based on the new gradually varied function to refine the $F$. This method will also work. The result will be smoothed to the order that we choose.

## 5. Experiments

In this section, we present three different experiments and applications. The first experiment uses the gradually varied function to fit the data for ground water distribution. The second compares the smoothness of reconstruction. At the end of the section, we show some examples for fitting continuous and smooth functions on manifolds.

**First experiment -- Ground water distribution**
Two sets of real data are tested. Each data set is of ground water distribution near Norfolk, Virginia. The first set consists of 10 sample data points of groundwater distribution and we call this the raw data. The format of the data is as follows:

```
Value    Latitude        Longititude
4.65     36.62074879     -76.10938540
8.60     36.65792848     -76.55772770
75.12    36.70764998     -76.12937859
208.26   36.68320624     -76.91329390
10.04    36.72371439     -76.02054720
…….
```

The gradually varied fitting result is shown in Fig. 5.1 and the second data set containing 29 sample points is shown in Fig. 5.2.

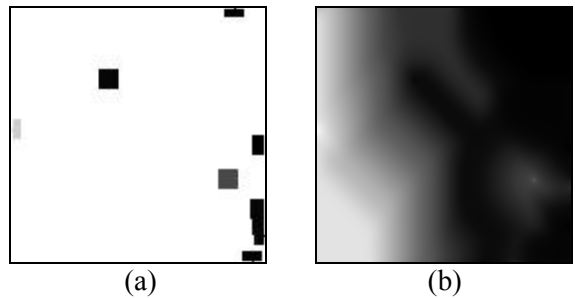

(a)          (b)

Fig. 5.1. Norfolk, VA Groundwater distribution calculated by gradually varied surfaces.  (a) Using10 sample points. (b) The fitted result.

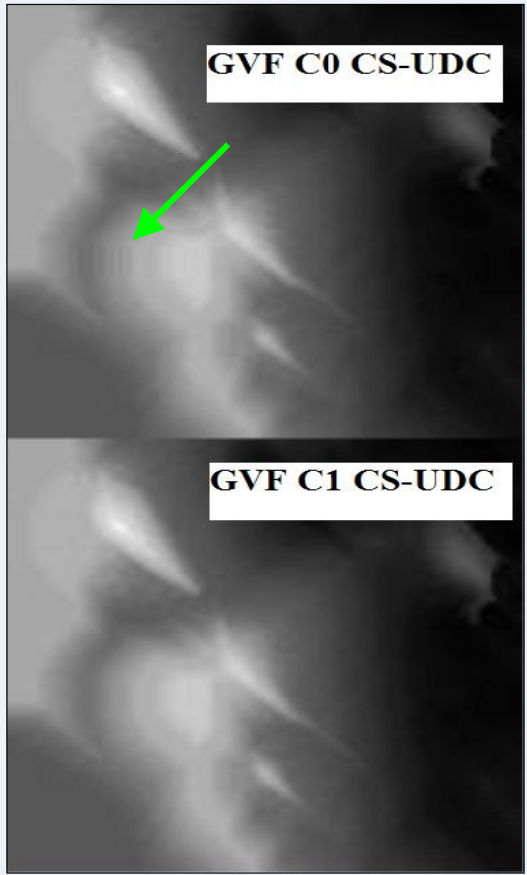

Fig. 5.2. The picture is the result of fitting based on 29 sample points. The first image is a "continuous" surface and the second is the "first derivative." The arrow indicates the interesting area that disappears in the second image, i.e. the vertical lines are removed.

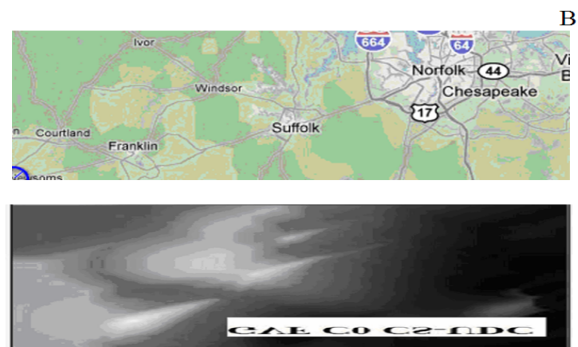

Dimensions (Latitude, Longitude)
A = (36.62074879, -77.17746768)
B = (36.92515020, -76.00948530)

Fig. 5.3. The map and ground water data

Fig. 5.3 shows a good match found between the ground water data and the region's geographical map. The brightness of the pixels indicates the water's depth from the surface. In mountainous areas, the groundwater level is lower in general. Some mismatches may be caused by not having enough sample data points (wells). The second image is done using the 29 points and is the same image as Fig. 5.2. However, this image is rescaled to match the geographical location.

**Second experiment -- Comparison of the smoothness of reconstruction**
This example uses the Taylor expansion formula to obtain the results. The original data is still the 29 points used above.

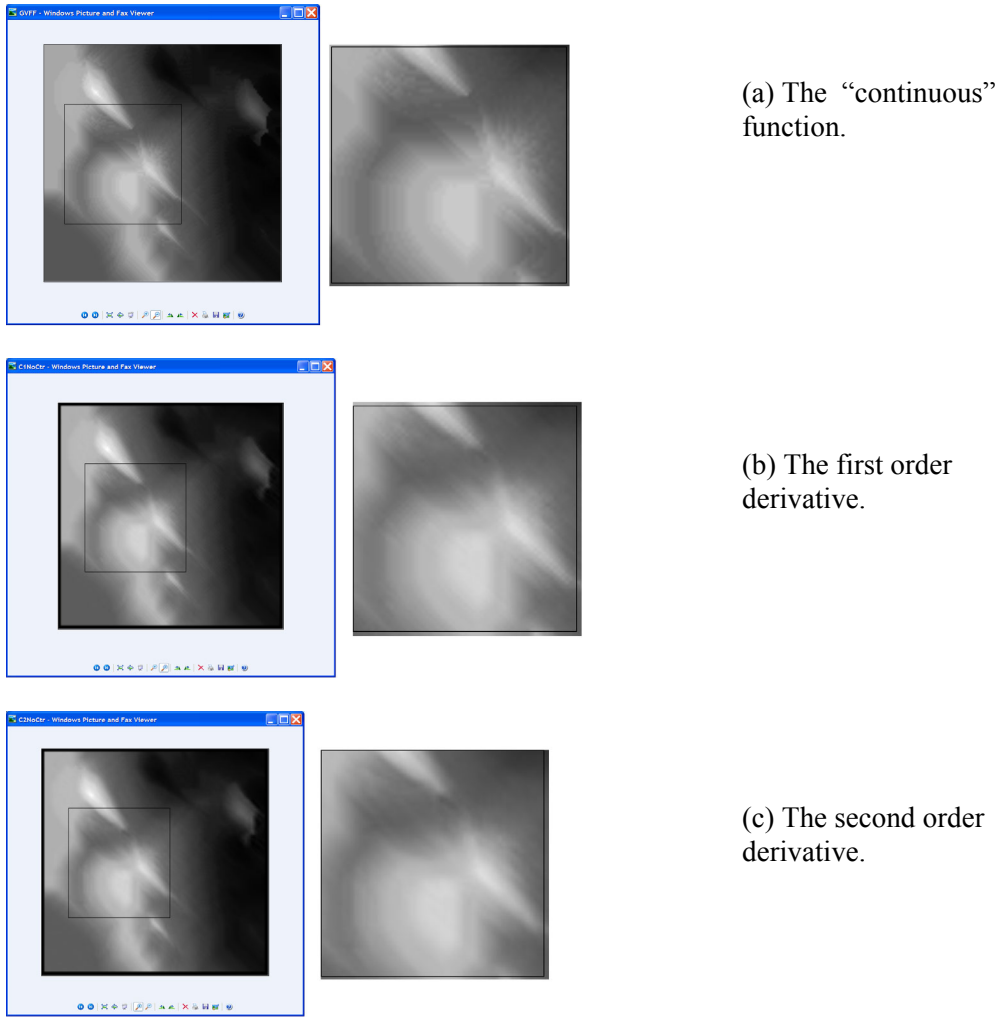

(a) The "continuous" function.

(b) The first order derivative.

(c) The second order derivative.

Fig. 5.4. Comparison for $C^{(0)}$, $C^{(1)}$ and $C^{(2)}$. The image on the right is an enlargement of the boxed region of the image to the left. One can see that the rightmost image is smoother in (c) than it is in (a).

**Third example -- Continuous and smooth functions on manifolds**
A gradually varied surface reconstruction does not rely on the shape of the domain and it is not restricted by simplicial decomposition. As long as the domain can be described as a graph, our algorithms will apply. However, the actual implementation will be much more difficult. In the above sections, we have discussed two types of algorithms for a rectangle domain. One is the complete gradually varied function (GVF) fitting and the other is the reconstruction of the best fit based on the gradual variation and finite difference methods.

The following is the implementation of the method for digital-discrete surface fitting on manifolds (triangulated representation of the domain). The data comes from a modified example in Princeton's 3D Benchmark data sets.

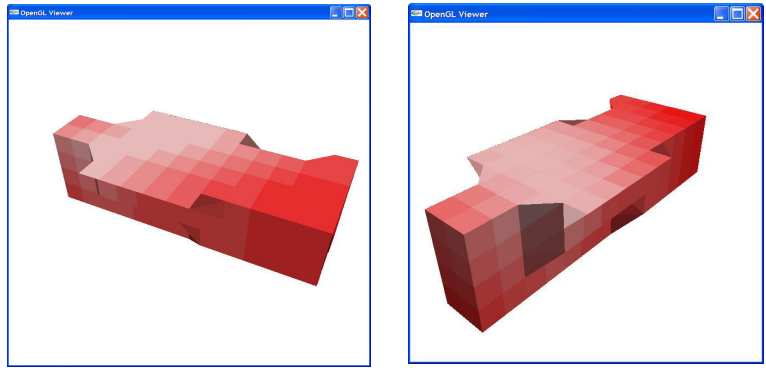

(a)                              (b)

Fig. 5.4. Gradually varied function on manifolds: (a) Fitting using seven points. (b) Harmonic fitting using (a).

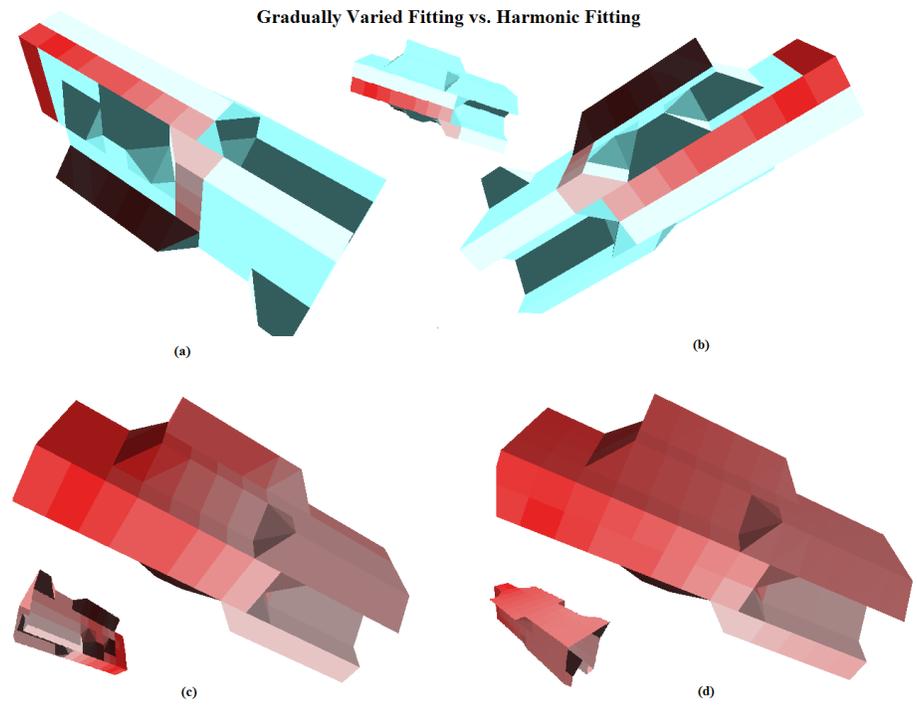

Fig. 5.5. More fitting examples in digital manifolds

We will have four algorithms related to continuous (and smooth) functions on manifolds. This is because we have 4 cases: (1) ManifoldIntGVF: The GVF extension on point space, corresponding to Delaunay triangulations; the values are integers. (2) ManifoldRealGVF: The GVF fitting on point space, the fitted data are real numbers. (3) ManifoldCellIntGVF: The GVF extension on face (2D-cell) space, corresponding to Voronoi decomposition; the values are integers. (4) ManifoldCellRealGVF: The GVF fitting on face (2D-cell) space, the fitted data are real numbers.

## 6. Summary

McShane-Whitney Theorem says that a Lipschitz function *f* on a subset *J* of a connected set *D* in a metric space can be extended to a Lipschitz function *F* on *D*. McShane gave a constructive proof for the existence of the extension in [17]. He constructed a minimal extension (INF) that is Lipschitz. It is easy for someone to construct a maximum extension (SUP). In [3], we use *F=(INF+SUP)/2* as the so-called *McShane-Whitney mid function.* The result is shown below. (See Fig. 6.1)

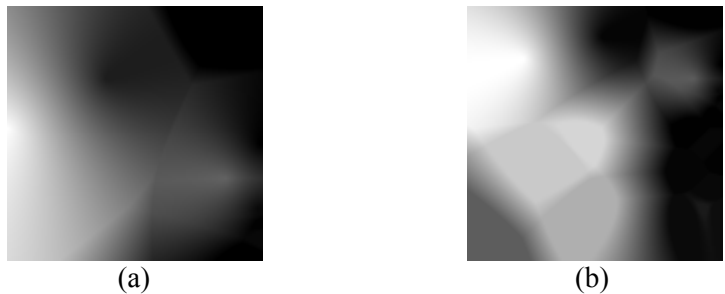

(a)  (b)

Fig. 6.1. McShane-Whitney mid extensions: (a) using the sample data set of Fig. 5.1, (b) using the sample data set of Fig. 5.2.

The fitting is dominated by the Lipschitz constant [3]. In this paper, we have shown the local Lipschitz function extensions. To get a smoothed function using gradual variation is a long time goal of our research. Some theoretical attempts have been made before, but struggled in the actual implementation. Fefferman et al have designed a refinement method [13].

The purpose of this paper is to present some actual examples and related results using the new algorithms we designed in [2]. The author welcomes other real data sets to further examine the new algorithms. The implementation code is written in C++. Li Chen's website can be found at www.udc.edu/prof/chen.

*Acknowledgements:* This research has been partially supported by the USGS Seed Grants through the UDC Water Resources Research Institute (WRRI) and Center for Discrete Mathematics and Theoretical Computer Science (DIMACS) at Rutgers University. Professor Feng Luo suggested the direction of the relationship between harmonic functions and gradually varied functions. Dr. Yong Liu provided many help in PDE. UDC undergraduate Travis Branham extracted the application data from the USGS database. Professor Thomas Funkhouser provided help on the 3D data sets and OpenGL display programs. The author would also like to thank Professor Charles Fefferman and Professor Nahum Zobin for their invitation to the Workshop on Whitney's Problem in 2009.

APPENDIX

Gradually Varied Definition:

Let function $f: D \rightarrow \{A_1, A_2,...,A_n\}$ and let $A_1 < A_2 < ... < A_n$. If $a$ and $b$ are adjacent in $D$, assume $f(a)=A_i$, then $f(b) = A_i$, $A_{i-1}$ or $A_{i+1}$. Point $(a,f(a))$ and $(b,f(b))$ are then said to be gradually varied.

A 2D function (surface) is said to be gradually varied if every adjacent pair is gradually varied.

Discrete Surface Fitting Definition:

Given $J \subset D$, and $f: J \rightarrow \{A_1, A_2,...A_n\}$, decide if there exists an F: $D \rightarrow \{A_1, A_2,...,A_n\}$ such that $F$ is gradually varied where $f(x)=F(x)$, and $x$ is in $J$.

An example using real numbers is:

Let D be a subset of real numbers $D \rightarrow \{1,2,.,.,.,n\}$. If $a$ and $b$ are adjacent in D such that $|f(a)-f(b)| \leq 1$ then point $[a,f(a)$ and $b(f(b)]$ is said to be gradually varied.

A 2D function (a surface) is said to be gradually varied if every adjacent pair are gradually varied.

There are three theorems necessary for this development:

**Theorem 2.1** (Chen, 1989) [8][9] The necessary and sufficient conditions for the existence of a gradually varied extension $F$ are: for all $x,y$ in $J$, $d(x,y) \geq |i-j|$, $f(x)=A_i$ and $f(y)=A_j$, where $d$ is the distance between $x$ and $y$ in $D$.

**Theorem 2.2** [7]: Any graph (or digital manifold) $D$ can normally immerse an arbitrary tree $T$.

**Theorem 2.3** [1]: For a reflexive graph $G$, the following are equivalent:
- (1) $G$ has the Extension Property.
- (2) $G$ is an absolute retract.
    - A space X is known as an *absolute retract* if for every normal space Y that embeds X as a closed subset, X is a retract of Y.
- (3) $G$ has the Helly property.

An alternative representation of the theorem is: For a discrete manifold $M$ the following are equivalent:
- (1) Any discrete manifold can normally immerse to $M$.
- (2) Reflexivized $M$ is an absolute retract.
- (3) $M$ has the Helly property.